\documentclass[12pt,oneside,final]{amsart}
\usepackage[cp1250]{inputenc}
\usepackage[T1]{fontenc}
\usepackage{amssymb}
\usepackage{color}
\renewcommand{\thefootnote}{}
\newtheorem{twr}{Twierdzenie}
\newtheorem{thr}[twr]{Theorem}
\newtheorem{prp}[twr]{Proposition}

\newtheorem{lm}[twr]{Lemma}

\newtheorem{crl}[twr]{Corollary}
\renewcommand{\thefootnote}{}

\newcommand{\po}{{\partial\Omega}}

\newcommand{\cn}{\mathbb{C}^n}

\title[$m$-subharmonic functions]{The smoothing of  $m$-subharmonic functions}
\author[S. Pli\'{s}]{Szymon Pli\'{s}}
\address{  Institute of Mathematics, Cracow University of Technology, Warszawska 24, 31-155
    Krak\'{o}w, Poland
}

\email{splis@pk.edu.pl}
\subjclass[2010]{ 35J66,  35B65,   32W50}
\keywords{complex Hessian equation, Richberg theorem, $m$-subharmonic function}

\begin{document}\thispagestyle{empty} \footnotetext{The author was partially supported by the NCN grant 2013/08/A/ST1/00312 }
\renewcommand{\thefootnote}{\arabic{footnote}}

\begin{abstract}
 We prove Richberg type theorem for $m$-subharmonic function. The main tool is the complex Hessian equation for which we obtain the existence of the unique smooth solution in strictly pseudoconvex domains.
\end{abstract}

\maketitle

\section{Introduction}

In the paper $M$ is a complex manifold with a fixed K\"{a}hler form $\beta$. Let $\Omega\subset M$ be a domain.  We say that a function $u\in\mathcal{C}^2(\Omega)$ is $m$-subharmonic if $(dd^cu)^k\wedge\beta^{n-m}\geq0$ for $k=1,\ldots,m$.   We say that a locally integrable function $$u:\Omega\rightarrow[-\infty,+\infty)$$ is $m$-subharmonic ($u\in\mathcal{SH}_m(\Omega)$) if $u$ is upper semicontinuous and $$dd^cu\wedge dd^cu_1\wedge\ldots\wedge dd^cu_{m-1}\wedge\beta^{n-m}\geq0,$$ for any  $\mathcal{C}^2$  $m$-subharmonic functions $u_1,\ldots u_{m-1}$. We say that $u\in\mathcal{SH}_m(\Omega)$ is strictly $m$-subharmonic if for any $\varphi\in\mathcal{C}^\infty_0(\Omega)$ there is $\varepsilon>0$ such that $u+\varepsilon\varphi\in\mathcal{SH}_m(\Omega)$. For K\"{a}hler form $\omega$ we  say that $u$ is (strictly) $\omega\;-\;m$-subharmonic if  $u+\rho$ is (strictly) $m$-subharmonic, where $\rho$ is a local potential of $\omega$. 

The definition of $m$-subharmonity was given in \cite{b} in the case  of $\beta$  the Euclidean K\"{a}hler form in $\cn$. In this case we can smooth $m$-subharmonic functions by the standard regularisation. In \cite{c} and \cite{d-k1}   definitions are given for any K\"{a}hler $\beta$, however in \cite{d-k1} authors assume formally stronger condition, that any $m$-subharmonic function is locally a limit of a decreasing sequence of smooth $m$-subharmonic functions. In this paper we prove a Richberg type theorem for $m$-subharmonic functions, which gives us  this two definitions coincide for continuous functions.

The main results are the following
\begin{thr}\label{globalnaaproksymacjaciaglych}
If $M$ is compact, $\omega$ is a K\"{a}hler form on $M$, $u\in\mathcal{C}(M)$, $u$ is $\omega\;-\;m$-subharmonic and $h>0$, then  there exists a strictly $\omega\;-\;m$-subharmonic function $\psi\in\mathcal{C}^\infty(M)$ such that $u\leq\psi\leq u+h$.
\end{thr}

\begin{thr}\label{aproksymacjaciaglych}
If $u,h\in\mathcal{C}(M)$, $h>0$ and $u$ is strictly $m$-subharmonic, then  there exists a strictly $m$-subharmonic function $\psi\in\mathcal{C}^\infty(M)$ such that $u\leq\psi\leq u+h$.
\end{thr}

We prove theorems above using methods from \cite{p2} and similarly as there the main tool is the existence of smooth solutions of the Dirichlet problem for the Hessian equation.
\begin{thr}\label{DP}
Let $\Omega\subset\cn$ be a (bounded) strictly pseudoconvex domain and let $dV$ be a volume form on some neighbourhood of $\bar\Omega$.	
 Then the Dirichlet problem
 \begin{equation}\label{eDP}
\left\{
\begin{array}{l}
    u\in \mathcal{SH}_m(\Omega)\cap \mathcal{C}^\infty(\bar\Omega)\\ 
    (dd^c u)^m\wedge\beta^{n-m}=dV \;\mbox{ in }\;\Omega\\
    u=\varphi\;\mbox{ on }\;\partial\Omega
\end{array}
\right.
\end{equation}
has an unique solution.
\end{thr}

For related results in case of $\beta=dd^c|z|^2$ see \cite{l,b}. Note also that the existence of smooth solution of the Hessian equation on a compact K\"{a}hler manifold was proved in \cite{d-k2}. 

Another application of Theorem \ref{DP} is the existence of a continuous solution of the homogeneous Hessian equation. As in \cite{b} one can prove the following
\begin{crl}
Let $\Omega\subset\cn$ be a (bounded) strictly pseudoconvex domain.	
 Then the Dirichlet problem
 \begin{equation}\label{eDP}
\left\{
\begin{array}{l}
    u\in \mathcal{SH}_m(\Omega)\cap \mathcal{C}(\bar\Omega)\\ 
    (dd^c u)^m\wedge\beta^{n-m}=0 \;\mbox{ in }\;\Omega\\
    u=\varphi\;\mbox{ on }\;\partial\Omega
\end{array}
\right.
\end{equation}
has an unique solution.
\end{crl}

 Recently, after the author had written this paper, Chinh and  Nguyen proved in \cite{c-n} that on a compact K\"{a}hler manifold any $\omega\;-\;m$-subharmonic function is a limit of a decreasing sequence of smooth $\omega\;-\;m$-subharmonic functions\footnote{For $\omega=\beta$ which is a standard assumption.}. In their proof they use solutions of the complex Hessian equation on Compact manifold instead of strictly pseudoconvex domains. However they use Theorem \ref{DP} to deal with $m$-extremal functions.

\section{Proof of Theorem \ref{DP}}
\label{ProofDP}

To prove Theorem \ref{DP} it is enough to obtain {\it a priori estimates} up to the second order. An uniform estimate  and a gradient estimate 
 are standard and the second order interior estimates follow from \cite{h-m-w}. The uniqueness follows from the comparison principle.

Our proofs of {\it a priori estimates} are rather standard and close to proofs in \cite{b} but in our situation we can not choose local coordinates such that vectors $\frac{\partial}{\partial z_1},\ldots,\frac{\partial}{\partial z_n}$ are perpendicular. Instead of this, (similarly as in \cite{p1}) we work with vector fields which are not necessary commutative.

In the proofs $\zeta_1,\ldots,\zeta_n$ is always a (local) orthonormal frame of $T^{1,0}$ i.e. $$\beta=2i\sum \zeta_p^\star\wedge \bar\zeta_p^\star,$$ where $\zeta_1^\star,\ldots, \zeta_n^\star,\bar\zeta_1^\star,\ldots, \bar\zeta_n^\star$ is a base of $(T_{\mathbb{C}}M)^\star$ dual to the base\\ $\zeta_1,\ldots,\zeta_n,\bar{\zeta_1},\ldots,\bar{\zeta_n}$ of $T_{\mathbb{C}}M$. Let us put for a smooth function $u$
$$u_{p\bar q}=\zeta_p\bar{\zeta_q}u=u_{\bar qp}+[\zeta_p,\bar{\zeta_q}]u$$
 and
$$A_{p\bar q}=A_{p\bar q}(u)=u_{p\bar q}-[\zeta_p,\bar{\zeta_q}]^{0,1}u,$$
where $X^{0,1}=\Pi^{0,1}(X)$. Then for a smooth function $u$ we have (see \cite{p}):
 $$dd^c u=2i\sum A_{p\bar q}\zeta_p^\star\wedge \bar\zeta_q^\star.$$

In this section we assume that $\Omega\Subset M$ is strictly pseudoconvex of class $\mathcal{C}^\infty$ with the defining function  $\rho$. All norms of functions are taken with respect to $\beta$ or more precisely with respect to a rimannian metric which is given by $g(X,Y)=-\beta(X,JY)$ for vector fields $X$, $Y$.

Now we recall some facts from \cite{b}.

  Let $f\in\mathcal{C}^\infty(\bar{\Omega})$ be such that $dV=f\beta^n$. Then locally our Hessian equation $(dd^c u)^m\wedge\beta^{n-m}=f\beta^n $ has a form:
$$ S_m(A_{p\bar q})=f,$$ where $S_m$ is the $m$-th  elementary symmetric function of eigenvalues of the matrix $(A_{p\bar q})$. For a matrix $B=(b_{p\bar q})$ we put $D_m(B)=(\frac{\partial S_m(B)}{\partial b_{p\bar q}})$ and we have
\begin{equation}\label{nierownoscG}tr(A_1D_m(A_2))\geq mS_m(A_1)^{1/m}S_m(A_2)^{(m-1)/m}, \end{equation}
\begin{equation}\label{wzornaf}tr(A_1D_m(A_1))= mS_m(A_1) \;,\end{equation}

for  $A_1,A_2$ such that $ S_k(A_i)\geq0 $ for $i=1,2$ and $k=1,\ldots,m$.

Put $(a_{p\bar q})=D_m(A_{p\bar q})$. Recall that the product of matrices $(a^{p\bar q})$ and $(A_{p\bar q})$ is a hermitian matrix, what implies that  for every $p$, $q$ \begin{equation}\label{dlaczego}a^{p\bar l}A_{q\bar l}=a^{k\bar q}A_{k\bar p}\;.\end{equation}
From (\ref{wzornaf}) we get 
 \begin{equation*}\label{pp}
    a^{p\bar q}XA_{p\bar q}=Xf\;.
 \end{equation*} We very often use the following elliptic operator
$$L=L_{\zeta}=a^{p\bar q}(\zeta_p\bar{\zeta_q}-[\zeta_p,\bar{\zeta_q}]^{0,1}).$$

 In the Lemmas we  specify exactly how  {\it a priori estimates} depend on $\rho$, $f$ and $\varphi$. We should  emphasize that they also depend strongly on   $\beta$. The notion $C(A)$ really means that $C$ depends on an upper bound for $A$. $C$ always depends on    $m(\rho)$ which is defined as the smallest constant $m>0$ such that $\beta\leq mi\partial\bar{\partial}\rho$ on $\Omega$.

In the proofs below $C$ is a constant under control, but it can change from a line to a next line.

Let us fix a point $P\in\po$. Now we  give the $\mathcal{C}^{1,1}$ estimate in a point $P$ (which not depends on $P$). We can assume that $P=0\in\cn$, $\beta(0)=dd^c|z|^2(0)$ and $\triangledown\rho=\frac{\partial}{\partial y_n}$. For $k=1\ldots,n$ we put $X_{2k-1}=\frac{\partial}{\partial x_k}$, $X_{2k}=\frac{\partial}{\partial y_k}$. The estimate of $XYu(P)$, where $X,Y$ are tangent to $\po$, follows from the gradient estimate.

\begin{lm}\label{druganabrzegu}
Let $X$ be a vector field on a neighborhood of $P$ tangent to $\po$ on $\po$. We have \begin{equation*}\label{druganabrzeg}
|X_{2n}Xu(P)|\leq C,
\end{equation*}
 where 
$C=C(\|\rho\|_{\mathcal{C}^{0,1}(\Omega)},\|f^{1/n}\|_{\mathcal{C}^{0,1}},\|\varphi\|_{\mathcal{C}^{2,1}(\Omega)} , \|X\|_{\mathcal{C}^{0,1}},\|u\|_{\mathcal{C}^{0,1}(\Omega)})$.
 \end{lm}

\textit{Proof:} Consider the function $$v=X(u-\varphi)+\sum_{k=1}^{2n-1}|X_k (u-\varphi)|^2-A|z|^2+B\rho.$$
Let $V\Subset U$ be a neighbourhood of $P$ and  $S=V\cap\Omega$. For $A$ large enough $v\leq0$ on $\partial S$.

Our goal is to show that for $B$ large enough we have $v\leq0$ on $\bar S$. Let $z_0\in S$ be a point where $v$ attains a maximum and let a frame $\zeta_1,\ldots,\zeta_n$ be orthonormal  such that $\zeta_1(z_0),\ldots,\zeta_{n-1}(z_0)\in\rm{Span}_\mathbb{C}(X_1,\ldots,X_{2n-2})$ and $A_{p\bar q}(z_0)=0$ for $p<q<n$. From now on all formulas are assumed to hold at $z_0$. It is clear that:
$$ \sum a^{p\bar p}\leq CL(\rho) $$ and
$$L(-X\varphi-A({\rm dist}(P,\cdot))^2)\geq-C\sum a^{p\bar p},$$ hence for $B$ large enough
$$L(B\rho-X\varphi-A({\rm dist}(P,\cdot))^2)\geq \frac{B}{2}\sum a^{p\bar p}.$$
To estimate $L(Xu+\sum_{k=1}^{n}|X_k (u-\varphi)|^2)$ let us first consider $Y\in\{X,X_1,\ldots, X_{n}\}$ and  calculate 
$$L(Yu)=a^{p\bar q}(\zeta_p\bar\zeta_q Yu-[\zeta_p,\bar\zeta_q]^{0,1}Yu)$$ $$=Y f +a^{p\bar q}(\zeta_p[\bar\zeta_q, Y]u+[\zeta_p, Y]\bar\zeta_qu-[ [\zeta_p,\bar\zeta_q]^{0,1},Y]u).$$
There are $\alpha_{q,k},\beta_{q,k}\in\mathbb{C}$ such that
$$[\bar\zeta_q, Y]=\sum_{k=1}^{n}\alpha_{q,k}\bar\zeta_k+\sum_{k=1}^{2n-1}\beta_{q,k}X_k$$
and so
$$a^{p\bar q}\zeta_p[\bar\zeta_q, Y]u=\sum_{k,l}\alpha_{k,l}(a^{p\bar k}A_{p\bar l})+\sum_{k=1}^{2n-1}a^{p\bar q}\beta_{q,k}\zeta_pX_ku+a^{p\bar q}Z_{pq} u,$$
 where $Z_{pq}$ are  vector fields under control. For $k<n$, by (\ref{dlaczego}) we get $$a^{p\bar k}A_{p\bar n}=a^{n\bar q}A_{k\bar q}=a^{n\bar k}A_{k\bar k}+a^{n\bar n}A_{n\bar k}$$ and by (\ref{wzornaf}) $$a^{p\bar n}A_{p\bar n}=f-\sum_{q<n,p}a^{p\bar q}A_{p\bar q}.$$
 This gives us
$$|a^{p\bar q}\zeta_p[\bar\zeta_q, Y]u|\leq C \sum_{k<2n,q}a^{p\bar q}(1+|\zeta_pX_ku|).$$
 In a similar way we can estimate $a^{p\bar q}[\zeta_p, Y]\bar\zeta_qu$ and we obtain
$$|L(Yu)|\leq C \sum_{k<2n,q}a^{p\bar q}(1+|\zeta_pX_ku|).$$
Therefore   we obtain
$$L(Xu+\sum_{k<2n}|X_k (u-\varphi)|^2)$$
   $$\geq a^{p\bar q}\sum_{k<2n}(\zeta_pX_k (u-\varphi))(\bar\zeta_{ q}X_k (u-\varphi))-C \sum_{k<2n,q}a^{p\bar q}(1+|\zeta_pX_ku|)$$
  $$\geq \sum_{k<2n}a^{p\bar q}\zeta_pX_k u\bar\zeta_qX_k u-C \sum_{k<2n,q}a^{p\bar q}(1+|\zeta_pX_ku|).$$

Now for $B$ large enough, since by the Schwarz inequality $L(v)(z_0)>0$,  we have  contradiction with maximality of $v$. Hence $v\leq0$ on $S$ and so $X_{2n}Xu(P)\leq C$ $\;\Box$.

\begin{lm}\label{druganormalnanabrzegu}
 We have \begin{equation}\label{drugannabrzegu}
\|\frac{\partial^2u}{\partial y_n^2}(0)\|\leq C,
\end{equation}
 where  
$$C=C(\|\rho\|_{\mathcal{C}^{2,1}(\Omega)},\|f^{1/n}\|_{\mathcal{C}^{0,1}},\|f^{-1}\|_{L^\infty({\Omega})},\|\varphi\|_{\mathcal{C}^{3,1}(\Omega)} , \|u\|_{\mathcal{C}^{0,1}(\Omega)}).$$
 \end{lm}

\textit{Proof:} Let a frame $\zeta_1,\ldots,\zeta_n$ be orthonormal  such that vectors  $\zeta_1,\ldots,\zeta_{n-1},X$ are tangent to $\partial\Omega$ on $\partial\Omega$ where $\zeta_n=X-iJX$. Let us put $M'=(m_{p\bar q})_{1\leq p,q\leq m-1}$  for a matrix  $M=(m_{p\bar q})$, $A'=(A_{p\bar q})'$ and $S_{m-1}'=S_{m-1}(A')$. We can write $$ f=A_{n\bar n}S'_{m-1}+O(1).$$
We may assume that $S'_{m-1}|_{\partial\Omega}$ has the minimum at $0$. By (\ref{wzornaf}) and (\ref{nierownoscG})
$$m tr\left(B_0(A'(z)-A'(0))\right)\geq S'(z)-S'(0)\geq0$$ for $z\in\partial\Omega$,
where $$B_0=D_{m-1}(A'(0)).$$
This gives us $$w=Nu(z)\psi(z)-Nu(0)\psi(0)+m tr\left(B_0(\Phi(z)-\Phi(0))\right)\geq0,$$ where $N=JX$ $$\psi=\frac{tr(B_0(\zeta_p\bar\zeta_q\rho-[\zeta_p,\zeta_q]^{0,1}\rho)')}{trB_0}$$
and $$\Phi=\frac{(\zeta_p\bar\zeta_q\varphi-[\zeta_p,\zeta_q]^{0,1}\varphi)'-(\zeta_p\bar\zeta_q\rho-[\zeta_p,\zeta_q]^{0,1}\rho)'N\varphi}{trB_0}.$$
Similarly as in the proof of the previous Lemma we can prove that choosing $A$, $B$ large enough (but under control) a barrier function $$v=-w+\sum_{k=1}^{2n-1}|X_k (u-\varphi)|^2-A|z|^2+B\rho$$ is non positive in $U\cap\Omega$, where $U$ is some neighbourhood of $0$. We thus obtain $\frac{\partial^2u}{\partial y_n^2}(0)\psi\leq0$ which gives (\ref{drugannabrzegu}). 
 $\;\Box$

\section{Approximation}
The following lemma generalizes lemma 3.7 and proposition 5.1 from \cite{lu1}.
\begin{lm}\label{lepkosc}
An uppersemicontinuous function $H$ is $m$-subharmonic iff for any $p\in\Omega$ and any $\mathcal{C}^2$ function $\varphi\geq H$ such that $\varphi(p)=H(p)$, we have $H_k(\varphi)(p)\geq 0$ for $k=1,\ldots,m$.
\end{lm}

{\it Proof:} Let $\beta_1,\ldots,\beta_{m-1}$ be smooth $m$-positive $(1,1)$ forms. A $(n-1,n-1)$-form $\Omega=\beta_1\wedge\ldots\wedge\beta_{n-1}\wedge\beta^{n-m}$ is a closed positive form  and there is a positive form $\omega$ such that $\Omega=\omega^{n-1}$. In local coordinates we have $\omega=i\sum_{p,q}g_{p\bar q}dz_p\wedge d\bar z_q$ for some hermitian matrix $(g_{p\bar q})\geq0$. It is easy to check that $dd^cH\wedge\Omega\geq0$ iff $g^{p\bar q}H_{p\bar q}\geq0$. Now from the theory of linear elliptic operators (see section 9 in \cite{h-l}) we obtain that  $dd^cH\wedge\Omega\geq0$ iff $dd^c\varphi\wedge\Omega(p)\geq0$ for any $\mathcal{C}^2$ function $\varphi\geq H$ such that $\varphi(p)=H(p)$. $\;\Box$

We need the following version of the comparison principle
\begin{prp}\label{cp'}
Suppose that $\Omega$ admits a bounded, smooth strictly plurisubharmonic function and   $u$, $v\in\mathcal{C}^2\cap\mathcal{SH}_m(\Omega)$ are such that  $H_m(u)\geq H_m(v)$. Then for any $H\in\mathcal{SH}_m(\Omega)$, an inequality $$\varlimsup_{z\rightarrow z_0}(u+H-v)\leq0$$ for any $z_0\in\po$ implies $u+H\leq v$ on $\Omega$.
\end{prp}

{\it Proof:} Let us assume that $H_m(u)> H_m(v)$ and a function $u+H-v$ attains a maximum in a point $p\in\Omega$. Using above Lemma (for $\varphi=v-u+A$ where $A$ is such that $\varphi(p)=H(p)$) we get that $H_k(v-u)(p)\geq0$ for $k=1,\ldots,m$. This gives $H_m(v)(p)\geq H_m(v-u)(p)+H_m(u)(p) \geq H_m(u)$ which is a contradiction.
 The general case ($H_mu\geq H_mv$) we obtain as usually from the case  above by adding  to $u$ a small, smooth, negative strictly $m$-subharmonic function. $\;\Box$

Exactly as in \cite{p2} (see proposition 3.3 there) we can prove the following
\begin{lm}\label{lokalnerozwiazaniePD}
Let  $u$ be a continuous strictly $m$-subharmonic function. If $ U\Subset M$ is a smooth strictly pseudoconvex domain 
 and  $K\Subset U$, then there is $v\in\mathcal{C}^\infty(\bar U)$  strictly $m$-subharmonic function on $U$ such that $v<u$ on $\partial U$ and $v>u$ on $K$.
\end{lm}

{\it Proof of Theorem \ref{globalnaaproksymacjaciaglych}:} We can assume that $u$ is a strictly $\omega\;-\;m$-subharmonic function. Let us consider two open finite coverings $\{U_k\}$, $\{U_k'\}$, $k=1,\ldots,N$ of $M$ such that for every $k$:\\ $\bullet$  a domain $U_k$ is smooth strictly pseudoconvex,\\ $\bullet$ $\bar U_k'\subset U_k$,\\ $\bullet$ there is a function $\rho_k$ in neighbourhood of $U_k$ with $dd^c\rho_k=\omega$ such that $\sup_{U_k}u+\rho_k<h+\inf_{U_k}u+\rho_k$.\\
By Lemma \ref{lokalnerozwiazaniePD} there are smooth strictly $m$-subharmonic  functions such that $v_k>u+\rho_k$ on $\bar U_k'$ and $v_k<u+\rho_k$ on $\partial U_k$. Then $v_k<u+\rho_k+h$. For any $k$ we can easily modify outside $\bar U_k'$ (and extend) a function $v_k-\rho_k$ to a function $u_k\in\mathcal{C}^\infty(M)$ such that:\\$\bullet$ $u_k<u+h$ on $M$, \\$\bullet$ $u_k>u$ on $\bar U_k'$,\\ $\bullet$ $u_k<u$ on $M\setminus U_k$ and \\ $\bullet$ $u_k $ is strictly $m$-subharmonic function on set $\{u_k-u>\frac{1}{2}\inf_{ U_k}u_k-u\}$.

Let $j\in\mathbb{N}$. Define $$\psi=\frac{1}{j}\log(e^{ju_1}+\ldots+e^{ju_N}). $$ Observe that $\psi>u$ and for $j$ large enough $\psi$ is a strictly $\omega\;-\;m$-subharmonic function with $\psi<u+h$. $\;\Box$

Using Proposition \ref{cp'} and Lemma \ref{lokalnerozwiazaniePD} we can prove Theorem \ref{aproksymacjaciaglych} in exactly the same way as Theorem 3.1 in \cite{p2}.


$\newline$\textbf{Acknowledgments.} The author would like to express his
gratitude to Z. B\l ocki, L. H. Chinh, S. Dinew, S. Ko\l odziej, N. C. Nguyen for helpful discussions.

\end{document}